\documentclass{article}
\usepackage{amssymb,amsmath,amsthm,verbatim,IEEEtrantools,macros,bm,authblk}
\usepackage[margin=1in]{geometry}
\DeclareMathOperator{\cyl}{cyl}
\DeclareMathOperator{\supp}{supp}
\DeclareMathOperator{\diam}{diam}
\DeclareMathOperator{\Exact}{Exact}
\newcommand{\ltyp}{\lambda_{\text{typ}}}
\newcommand{\ltypn}{\lambda_{\text{typ},n}}
\newcommand{\lexc}{\lambda_{\text{exc}}}
\newcommand{\lexcn}{\lambda_{\text{exc},n}}

\newcommand{\dimf}{\text{dim}_F}
\title{Fourier decay of measures supported on sets of numbers with consecutive partial quotients belonging to a given set}
\author{Robert Fraser}
\affil{Wichita State University}
\begin{document}
\maketitle
\begin{abstract}
We consider measures supported on sets of irrational numbers possessing many consecutive partial quotients satisfying a condition based on the previous partial quotients. We show that under mild assumptions, such sets will always support measures whose Fourier transform decays to zero.
\end{abstract}
\section{Introduction}
\subsection{The theory of metric Diophantine approximation}
The theory of metric Diophantine approximation concerns metrical properties of sets arising in Diophantine approximation. For example, let $\psi$ be a decreasing function and consider the \textbf{well-approximable numbers} $E(\psi)$ defined by
\[\left\{\left|x - \frac{p}{q} \right| \leq \psi(q) \quad \text{for infinitely many rational numbers $\frac{p}{q}$}\right\}.\]
The study of the metrical properties of $E(\psi)$ goes back to an old result of Khintchine who gives a condition for which $E(\psi)$ has zero Lebesgue measure; a proof is provided by e.g. Bernik and Dodson \cite{BernikDodson99}. When $\psi(q) = q^{-\tau}$, Jarník \cite{Jarnik29} and Besicovitch \cite{Besicovitch34} compute the Hausdorff dimension of $E(\psi)$. Dodson \cite{Dodson92} computes the Hausdorff dimension of $E(\psi)$ for a more general family of decreasing functions $\psi$ in terms of the lower order of $\psi$ at infinity. Many results of this type are encapsulated by the mass transference principle of Beresnevich and Velani \cite{BeresnevichVelani06A}.

The theory of Diophantine approximation is intricately tied to the notion of partial fractions. If $x$ is an irrational number with continued fraction expansion $[a_0; a_1, a_2, \ldots]$, then the convergents $\frac{p_i}{q_i} =  [a_0; a_1, a_2, \ldots, a_i]$ are the best approximants to $x$ in the sense that $\norm{q_i x} < \norm{qx}$ for all $q < q_i$, where $\norm{\cdot}$ denotes the distance to the nearest integer. If the sequence $\{a_i\}_{i=0}^{\infty}$ is bounded by some number $N$, then the number $x$ is said to be \textbf{badly-approximable}. If $\mathcal{A} \subset \mathbb{N}$ is a finite set with at least two elements, we write $B(\mathcal{A})$ for the set of badly approximable numbers such that each partial quotient is contained in $\mathcal{A}$. The Hausdorff dimension of $B(\mathcal{A})$ is computed by Good \cite{Good41}.

A set of numbers in Diophantine approximation that combines some of the properties of the well-approximable numbers and some properties of the badly-approximable numbers is the set of \textbf{numbers approximable to exact order} $\Exact(\psi)$ introduced by Bugeaud \cite{Bugeaud03}. The number $x$ is said to belong to $\Exact(\psi)$ if the following two conditions hold:
\begin{enumerate}[(a)]
\item $\left| x - \frac{p}{q} \right| \leq \psi(q)$ for infinitely many rational numbers $\frac{p}{q}$, and
\item $\left| x - \frac{p}{q} \right| \geq (1 - o(1)) \psi(q)$ for all rational numbers $\frac{p}{q}$.
\end{enumerate}
Bugeaud shows that if $\psi$ is a function such that $q^2 \psi(q)$ is nonincreasing and such that $\sum_q q \psi(q)$ converges, then $\Exact(\psi)$ and $E(\psi)$ have the same Hausdorff dimension. 
\subsection{Fourier analysis and Diophantine approximation}
It is interesting to consider the Fourier transform of measures supported on sets arising from Diophantine approximation. An energy integral version of Frostman's lemma implies that the Hausdorff dimension of a Borel set $E$, denoted $\dim E$, is the supremum of those values of $s$ such that $E$ supports a Borel probability measure $\mu$ with
\begin{equation}\label{Frostman}
\int |\widehat{\mu}(\xi)|^2 |\xi|^{s - n} < \infty,
\end{equation}
so the Hausdorff dimension of a Borel set $E$ governs the $L^2$-averaged decay of measures supported on $E$. 

In contrast, the \textbf{Fourier dimension} of a Borel set $E$, denoted $\dimf E$, is the supremum of those values of $s$ such that $E$ supports a Borel probability measure $\mu$ with
\begin{equation}\label{Fourierdimension}
\widehat{\mu}(\xi) \lesssim (1 + |\xi|)^{-s/2}.
\end{equation}
Because the condition \eqref{Fourierdimension} implies the condition \eqref{Frostman}, it follows that $\dimf E \leq \dim E$ for any Borel set $E$. The value of $\dimf E$ can be anything from $0$ to $\dim E$ as has been observed by Körner in dimension $1$ \cite{Korner11}.

The Fourier dimension of the set of well-approximable numbers has been computed by Kaufman \cite{Kaufman81}, and the Fourier dimension of the badly-approximable numbers $B(\mathcal{A})$ has been shown to be positive for any set $\mathcal{A}$ with at least two elements by Sahlsten and Stevens \cite{SahlstenStevens24}; previous results in this direction are due to Kaufman \cite{Kaufman80}, Queffélec and Ramaré \cite{QueffelecRamare03}, and Hochman and Shmerkin \cite{HochmanShmerkin15}. If $\psi(q) = q^{-\tau}$, then the set $\psi$-exact order approximable numbers have been shown to have positive Fourier dimension by the author and Wheeler \cite{FraserWheeler22} \cite{FraserWheeler24}.

A weaker condition than \eqref{Fourierdimension} is given by Rajchman measures. A Borel probability measure $\mu$ is called a \textbf{Rajchman measure} if 
\[|\widehat{\mu}(\xi)| \to 0 \text{ as $|\xi| \to \infty$}.\]

Bluhm \cite{Bluhm00} shows that the set of Liouville numbers supports Rajchman measures. This result is strengthened by Polasek and Réla \cite{PolasekRela24}, who obtain explicit decay conditions for the Fourier transform of measures supported on the set of Liouville numbers, and the author and Nguyen \cite{FraserNguyen24}, who show that $E(\psi)$ supports Rajchman measures for any decreasing function $\psi$.

Rajchman measures and Fourier dimension are important in number theory for their connection to normal numbers. Recall that a number $x$ is called \textbf{normal} if the fractional part of $a^j x$ is equidistributed modulo $1$ for any integer $a \geq 2$. A result of Davenport, Erdős, and Leveque \cite{DavenportErdosLeveque63} implies that if $\mu$ is a measure and 
\begin{equation}\label{normalcond}
|\widehat{\mu}(\xi)| \lesssim (\log \log (10 + |\xi|))^{-(1 + \epsilon)}
\end{equation}
for some $\epsilon > 0$, then the support of $\mu$ must contain normal numbers. A proof of this result can be found in the article of Pollington, Velani, Zafeiropoulos, and Zorin \cite{PollingtonVelaniZafeiropoulosZorin22}.

\subsection{Techniques for constructing measures with rapid Fourier decay on sets in Diophantine approximation}
Broadly, there are two strategies for constructing Rajchman measures on sets of numbers in Diophantine approximation. One strategy that is particularly well-suited to the well-approximable numbers is to construct a measure whose support is concentrated near a large number of arithmetic progressions at many different scales. This strategy is used in constructions of measures with rapid Fourier decay supported on $E(\psi)$ such as that of Kaufman \cite{Kaufman81}, Bluhm \cite{Bluhm00}, Hambrook \cite{Hambrook14} \cite{Hambrook19}, the author and Hambrook \cite{FraserHambrook20} \cite{FraserHambrook23}, the author, Hambrook, and Ryou \cite{FraserHambrookRyou24}, and the author and Nguyen \cite{FraserNguyen24}. This technique is also used by the author and Wheeler \cite{FraserWheeler22} to estimate the Fourier dimension of $\Exact(\psi)$ in certain cases.  
The second main strategy, which is more adapted to the badly-approximable numbers, involves constructing a measure using properties of the continued fraction expansion. This technique is employed to construct measures supported on the badly approximable numbers with rapid Fourier decay by Kaufman \cite{Kaufman80} and Queffélec and Ramaré \cite{QueffelecRamare03}. These authors construct a periodic probability measure on bounded integer sequences that pushes forward to a measure supported on the real numbers. The periodicity of the measure guarantees a self-similarity in the pushforward that allows the authors to estimate the Fourier transform of the measure. A variant of this technique is employed by Hochman and Shmerkin \cite{HochmanShmerkin15}, Jordan and Sahlsten \cite{JordanSahlsten16} and Sahlsten and Stevens \cite{SahlstenStevens24}, who use the dynamical properties of the Gauss map to construct measures with desired properties.

We focus our attention on a variant of Kaufman's argument introduced by the author and Wheeler \cite{FraserWheeler24}. Rather than considering a purely periodic measure on bounded integer sequences, we introduce a set of sequences in which the majority of entries are bounded by some number $N$, but with a sparse set of exceptional partial quotients that are much larger than $N$. Although the pushforward of this measure is not self-similar, the failure of the measure to be self-similar can be carefully controlled. This technique is further used by Tan and Zhou \cite{TanZhou24} to estimate the Fourier dimension of sets of numbers $x$ in which there are infinitely many consecutive pairs of partial quotients $(a_n, a_{n+1})$ in which the product $a_n a_{n+1}$ is larger than some quantity $\Phi$, but such that each $a_n$ is bounded above by $\Phi$. Thus Tan and Zhou consider a set whose definition involves not only the behavior of individual partial quotients, but also the behavior of \textit{consecutive pairs} of partial quotients. This work will deal with a set defined by the behavior of arbitrarily long tuples of partial quotients.

\subsection{Statement of main result}
We discuss sets of numbers whose partial quotients satisfy a very general type of condition. We show under mild assumptions that sets of such numbers must support Rajchman measures, and obtain a quantitative estimate for the decay of measures supported on such sets.

Let $\mathbb{N}^*$ be the set of finite sequences of positive integers. Let $S : \mathbb{N}^* \to 2^{\mathbb{N}}$ be such that $S(a_0, \ldots, a_i) \neq \emptyset$ for any finite sequence $(a_0, \ldots, a_i) \in \mathbb{N}^*$. We will refer to such a map as an \textbf{assignment of sets to partial quotients}.

\begin{mydef}\label{eskdef}
Given such a mapping $S$, we define $E(S,k)$ to be the set of $x = [a_0; a_1, \ldots]$ such that there exist infinitely many $k$-tuples $a_{i_n}, \ldots, a_{i_n + k - 1}$ of consecutive partial quotients such that $a_i \in S(a_0, \ldots, a_{i-1})$ for $i = i_n, \ldots, i_n + k  - 1$.
We define the set $E(S, \infty) = \bigcap_{k=1}^{\infty} E(S,k)$. 
\end{mydef}

We highlight the example of numbers $x$ that are well-approximated by arbitrarily long sequences of consecutive convergents in the continued fraction expansion.
\begin{myex}
Let $\psi \geq 0$ be a decreasing function such that $q^2 \psi(q)$ has limit zero. We will define $E(\psi, k)$ to be the set of points $x$ such that there exist infinitely many $k$-tuples of \textbf{consecutive} convergents $\{\frac{p_{i,j}}{q_{i,j}}\}_{1 \leq i < \infty, 1 \leq j \leq k}$ such that 
\[\left|x - \frac{p_{i,j}}{q_{i,j}} \right| \leq \psi(q_{i,j}).\]
We write $E(\psi, \infty)$ for $\bigcap_{j = 1}^{\infty} E(\psi, k)$. 

Let $q(a_0, \ldots, a_n)$ denote the denominator of $[a_0; a_1, \ldots, a_n]$. If we define the mapping $S(a_0, \ldots, a_{n})$ by $S(a_0, \ldots, a_{n}) = \mathbb{N} \cap [q^2 \psi(q), \infty)$ then $E(\psi, k) \subset E(S,k)$.
\end{myex}
Of course, it is possible to place more exotic conditions on the mapping $S$. For example, $S(a_0, \ldots, a_n)$ could consist of the singleton set $\left\{\sum_{j=0}^{\infty} a_j \right\}$. In this case, the set $E(S,1)$ would consist of numbers that have infinitely many partial quotients equal to the sum of all of the previous partial quotients, and $E(S,\infty)$ consists of numbers that have arbitrarily long sequences of partial quotients, each of which is the sum of all of the previous partial quotients. 

\begin{mythm}
Let $S$ be any assignment of sets to partial quotients. Then $E(S, \infty)$ supports a Rajchman measure.
\end{mythm}

The fundamental strategy for this proof is similar to the strategy employed by the author and Wheeler \cite{FraserWheeler24} for estimating the Fourier dimension of the numbers approximable to exact order. We will construct a measure $\lambda$ on integer sequences that will push forward to a measure on $E(S, \infty)$ under the continued fraction map. For a sequence in the support of $\lambda$, the majority of the entries of the sequence will be integers bounded above by some number $N$; such entries will be called \textbf{typical}. However, there will be a small number of \textbf{exceptional} entries that might take a significantly larger value.

In the article \cite{FraserWheeler24} and the subsequent work by Tan and Zhou \cite{TanZhou24}, the strategy for controlling the exceptional partial quotients was to partition the measure into a small number of equivalence classes at a given scale, each of which could be treated in a similar way. This work will employ a simpler strategy: we will arrange our construction so that at any given scale $s$, only a small fraction of the measure will contain exceptional entries close to $s$. The piece of the measure corresponding to the exceptional partial quotients will be estimated trivially using the triangle inequality; the piece of the measure corresponding to a typical partial quotient will be estimated using Kaufman's argument.

As the trivial estimate for the exceptional portion of the measure will be quite large, we are unable to reproduce the Fourier dimension results of the previous result of the author with Wheeler \cite{FraserWheeler24} or the results of Tan and Zhou \cite{TanZhou24}. On the other hand, the flexibility of this method means that it can be used to solve a number of problems for which the previous sharper but more rigid methods do not apply.
\section{Notation and preliminaries}
\subsection{Notation}
Throughout the proof, we will fix a large integer $p$. A boldface letter such as $\bm{a}$ will denote a $p$-tuple of integers, and a letter with a ``vector arrow" such as $\vec a$ will denote a finite tuple of arbitrary length.

We let $\mathbb{N}^*$ denote the set of finite sequences of integers. The set $\mathbb{N}^{\infty}$ will denote the set of infinite sequences of integers. 

If $\vec a \in \mathbb{N}^*$, we write $\cyl(\vec a)$ for the set of elements of $\mathbb{N}^{\infty}$ beginning with $\vec a$ for some $\vec a \in X$. If $X \subset \mathbb{N}^*$, we write $X^*$ for the set of elements of $\mathbb{N}^*$ beginning with $\vec a$. We will write $X^{**}$ for the set of infinite sequences of integers beginning with $\vec a$ for some $\vec a \in X$. The space $\mathbb{N}^{\infty}$ will be viewed as a measurable space equipped with the cylinder $\sigma$-algebra generated by sets of the form $X^{**}$ where $X \subset \mathbb{N}^*$ is finite. Because the sets of the form $X^{**}$ form an algebra of sets, the Carathéodory extension theorem implies that any countably additive function defined for such sets extends to a unique measure on $\mathbb{N}^{\infty}$.  Therefore, in constructing a measure $\lambda$, it is sufficient to specify the value of $\lambda$ on such cylinder sets.

The function $g : \mathbb{N}^* \cup \mathbb{N}^{\infty} \to \mathbb{R}^+$ will denote the continued fraction map. That is, if $x = (a_0, a_1, \ldots, a_n)$ is a finite sequence, then $g(x)$ will denote the finite continued fraction $[a_0; a_1, \ldots, a_n]$ given by
\[[a_0; a_1, \ldots, a_n] = a_0 + \frac{1}{a_1 + \frac{1}{\cdots + \frac{1}{a_n}}},\]
and if $x = (a_0, a_1, \ldots)$ is an infinite sequence, then $g(x)$ will denote the infinite continued fraction $[a_0; a_1, a_2, \ldots]$ defined by
\[[a_0; a_1, a_2, \ldots] = a_0 + \frac{1}{a_1 + \frac{1}{a_2 + \frac{1}{\cdots}}}.\]

If $\lambda$ is any measure on the measurable space $\mathbb{N}^{\infty}$, we will write $g^{\#}\lambda$ for the pushforward of $\lambda$ under the continued fraction map $g$. Note that $g$ is a measurable function from $\mathbb{N}^{\infty}$ into $\mathbb{R}$ with the Borel $\sigma$-algebra, and a bijection from $\mathbb{N}^{\infty}$ into the set of irrational numbers.

For a finite sequence $(a_0, a_1, \ldots, a_i) =: \vec a \in \mathbb{N}^*$, we write $q(\vec a)$ for the least possible denominator of the continued fraction $g(\vec a)$, and we will write $q'(\vec a)$ for the least possible denominator of $(a_0, a_1, \ldots, a_{i-1})$. The denominators $q(\vec a)$ satisfy the recurrence
\[q(a_0, \ldots, a_i) = a_i q(a_0, \ldots, a_{i-1}) + q(a_0, \ldots, a_{i-2}).\]

Throughout this article, we write $A \lesssim B$ to indicate that $A$ is bounded above by a constant times $B$. This constant may depend on parameters such as $N$ and $\sigma$ appearing in the construction. We write $A \sim B$ for $A \lesssim B$ and $B \lesssim A$.
\subsection{Review of Kaufman's construction}
We will review Kaufman's construction \cite{Kaufman80} of a measure supported on the set of badly approximable numbers with polynomial Fourier decay. Kaufman's construction will be the basis for the construction in this work.

Kaufman constructs a measure $\lambda_K$ on the space supported on the space of infinite sequences of integers in $\{1, \ldots, N\}$. This measure is an infinite product of measures $\nu \times \nu \times \nu \times \cdots$ supported on $p$-tuples. 

Kaufman cites the following lemma of Rogers \cite{Rogers64}, whose statement we modify slightly. 
\begin{mylem}\label{joining}
Let $N \geq 2$ be an integer and suppose $(a_0, \ldots, a_j)$ and $(b_0, \ldots, b_k)$ be finite sequences with $1 \leq b_0 \leq N$. Then 
\[\left| \log q(a_0, \ldots, a_j, b_0, \ldots, b_k) - (\log q(a_0, \ldots, a_j) + \log q(b_0, \ldots, b_k)) \right| \leq C_N \]
where $C_N$ is a constant depending only on $N$.
\end{mylem}

Kaufman's measure $\nu$ (technically a slight variant) has the property that for any $(a_0, \ldots, a_{p-1}) \in \supp \nu$, we have
\[|\log q(a_0, \ldots, a_{p - 1}) - \sigma | \leq \frac{1}{10000} \sigma \]
where $\sigma$ is a number satisfying 
\[C_N \leq \frac{1}{10000} \sigma.\]
These two properties imply that for any $\vec a = (\bm{a}_0, \ldots, \bm{a}_i)$, where $\bm{a}_j \in \supp \nu$ for each $j$, that
\begin{equation}\label{qestimate}
\left|\log q(\vec a) - i \sigma \right| \leq \frac{1}{5000} i \sigma.
\end{equation}
Moreover, if $N$ is chosen sufficiently large, then $\nu$ has the property that if $(a_0, \ldots, a_{p-1}) \in \supp \nu$, then 
\begin{equation}\label{nuestimate}
\log \nu\{(a_0, \ldots, a_{p-1})\} \leq - \frac{198}{100} \sigma.
\end{equation}

Let $\nu^j$ denote the $j$-fold product $\nu \times \cdots \times \nu$. Because each element of $\supp \nu^j$ has $\nu^j$-measure no more than $2^{-j}$,  it follows from the pigeonhole principle that there exists a subset $T_j \subset \supp \nu \times \cdots \times \nu$ with $\left|\nu^j(T_j) - \frac{1}{2} \right| \leq 2^{-j}$. We will refer to $T_j$ as the \textbf{top half} and $T_j^c$ as the \textbf{bottom half} of $\supp \nu^j$. There is no particular significance to which elements of $\supp \nu^j$ are assigned to $T_j$ and which elements of $\supp \nu^j$ are assigned to $T_j^c$; we just need a convenient way to split $\supp \nu^j$ into two subsets of similar measure.

Combining  \eqref{qestimate} and \eqref{nuestimate}, we conclude that for $\vec a = (\bm{a}_0, \ldots, \bm{a}_{i-1})$ that
\begin{equation}\label{qnumeasure}
\log \nu^i(\{(\vec a)\}) \leq - \frac{196}{100} \log q(\vec{a}).
\end{equation}

Let $g^{\sharp} \lambda_K$ denote the pushforward of $\lambda_K$ under the continued fraction map. Let $B$ be any ball in $\mathbb{R}$. Because $N$ is a finite number, there exists $\vec a$ such that $B \subset g(\vec a)$ and such that the Lebesgue measure of $g(\vec a)$ is no more than an $N$-dependent constant times the diameter of $B$. Because the Lebesgue measure of $g(\vec a)$ is comparable to $q (\vec a)^{-2}$, we conclude from \eqref{qnumeasure} that for any ball $B$, we have the estimate
\begin{equation}\label{lambdafrostman}
g^{\sharp} \lambda_K(B) \lesssim \diam(B)^{98/100}.
\end{equation}
\section{Details of the construction}
We will let $\lambda_K, \nu, C_N$, and $\sigma$ be as in Kaufman's construction from the previous section. Recall that given a collection $X$ of finite sequences, we will write $X^*$ for the collection of all finite sequences beginning with an element of $X$.

We define a function $\phi : \mathbb{N}^* \to \mathbb{N}$ by 
\[\rho(\vec a) = \min S(\vec a),\]
we write $\rho_r$ for the map
\[\rho_r(\vec a) = \rho(\vec a, \rho(\vec a), \rho^2(\vec a), \ldots, \rho^{r-1}(\vec a)),\]
we write $\phi$ for the map
\[\phi(\vec a) = q((\vec a, \rho(\vec a))),\]
and we write $\phi_r$ for the map
\[\phi_r(\vec a) = q((\vec a, \rho(\vec a), \rho_2(\vec a), \ldots, \rho_{r}(\vec a)))\]
\begin{myex}
In the case of the $\psi$-well approximable numbers, recall that we chose $S(\vec a) = [\lceil q(\vec a)^2 \psi(q(\vec a)) \rceil, \infty)$. Suppose $\psi$ is a function such that $q^2 \psi(q)$ decreases to zero. In this case, the function $\rho(\vec a)$ is approximately
\[\rho(\vec a) = \frac{1}{q(\vec a)^2 \psi(q(\vec a))} + O(1)\]
and 
\[\phi(\vec a) =  q(\vec a) \rho(\vec a) + q'(\vec a),\]
where $q'(\vec a) < q(\vec a)$ is the previous denominator. Because $\frac{1}{q^2 \psi(q)} \to \infty$, we have 
\[\phi(\vec a) = (1 + o(1) ) \frac{1}{q(\vec a) \psi(q(\vec a))}.\]
The same argument shows that we have the recurrence

\[\rho_r(\vec a) = \frac{1}{\phi_{r-1}(\vec a)^2 \psi(\phi_{r - 1}(\vec a))} + O(1)\]
and
\[\phi_r(\vec a) = (1 + o(1))\frac{1}{\phi_{r-1}(\vec a) \psi(\phi_{r-1}(\vec a))}.\] 
Hence, we have the inequality
\[\phi_r(\vec a) < \frac{2}{\phi_{r-1}(\vec a) \psi(\phi_{r-1}(\vec a))}.\]
Writing $\Phi(q)$ for the function $\frac{1}{q \psi (q)}$, a simple induction shows that
\begin{equation}\label{phirestimate}
\phi_r(\vec a) \leq \Phi^{r}(2^{r-1} \vec q(a))
\end{equation}
where $\Phi^r$ denotes the $r$-fold composition of $\Phi$.
\end{myex}

To each $n \geq 1$, we will assign a weight $w_n = 2^{- \lfloor \log_2 n \rfloor}$. Observe that $\frac{1}{2n} \leq w_n \leq \frac{1}{n}$ for all $n$. Our measure $\lambda$ will be defined so that the partial quotient $a_{i_n}$ is exceptional for approximately a $w_n$-fraction of the measure.

Fix a nondecreasing sequence $\{r_n\}_{n=1}^{\infty}$ with $r_n \to \infty$ as $n \to \infty$. In practice, it is best to choose $r_n$ to be slowly growing, but this assumption is not necessary for the construction to work. We will inductively choose a sequence of \textbf{exceptional indices} $\{i_n\}_{n=1}^{\infty}$ and \textbf{exceptional locations} $\{X_n\}_{n=1}^{\infty}$ consisting of finite integer sequences. In our argument, $i_1$ will be chosen to be a very large index and $X_1 = \{\varepsilon\}$, where $\varepsilon$ denotes the empty string. At step $n$ for $n \geq 1$, we choose $i_{n+1}, X_{2n}$, and $X_{2n + 1}$. We will select the indices $\{i_n\}_{n=1}^{\infty}$ and the sets $\{X_n\}_{n=1}^{\infty}$ so that
\begin{enumerate}[(A)]
\item The $i_n$ are superlacunary in the sense that for any ratio $R > 0$, there exists $n_0(R)$ such that $\frac{i_{n+1}}{i_n + r_n} \geq R$ for $n \geq n_0(R)$.
\item The $i_n$ have the property that for every $n$:
\[i_{n+2} - i_{n+1} \geq \frac{100}{\sigma} \log \max_{(\bm{a}_0, \ldots, \bm{a}_{i_n}) \in X_n^*}\phi_{p r_n}((\bm{a}_0, \ldots, \bm{a}_{i_n})).\]
\end{enumerate}

The $X_n$ satisfy a few key properties:
\begin{enumerate}[(A)]
\setcounter{enumi}{2}
\item The sets $X_{2^t}^*, \ldots, X_{2^{t + 1} - 1}^*$ are disjoint for any $t \geq 0$. Moreover, for any $n$, the sets $X_{2n}^*$ and $X_{2n + 1}^*$ are contained in $X_n^*$.
\item The mass of each of $X_{2^t}^{**}, \ldots, X_{2^{t + 1} - 1}^{**}$ with respect to the measure $\lambda$ is $2^{-t}(1 + O(1))$. In other words, the mass of $X_n^{**}$ is $w_n(1 + O(2^{-i_1}))$ for every $n \geq 1$. Moreover, $\lambda(X_{2^t}^{**} \cup \cdots \cup X_{2^{t + 1} - 1}^{**}) = 1$, and for each $n$, $\lambda(X_n^{**}) = \lambda(X_{2n}^{**}) + \lambda(X_{2n + 1}^{**})$.
\end{enumerate}
We will describe the construction in detail now. The measure $\lambda$ will be defined on the space of infinite sequences of integers. Let $\mu$, $\nu$, etc. be as in the Kaufman construction. In Kaufman's construction and in ours, it is more convenient to work with $p$-tuples of entries rather than with individual entries. An infinite sequence of integers will therefore be written as $(\bm{a}_0, \bm{a}_1, \ldots)$ where each $\bm{a}_i$ is a $p$-tuple of integers.

The measure $\lambda$ will be defined via a mass-distribution procedure. Given a cylinder set $\cyl(\vec a)$ where $\vec a = (\bm{a}_0, \ldots, \bm{a}_{i-1})$, we will define two different ways of distributing the mass to cylinder sets $\cyl(\bm{a}_0, \ldots, \bm{a}_{i-1}, \bm{a}_i)$ at the next stage.

We say that $\bm{a}_i$ will be \textbf{chosen typically} if the mass is according to the measure $\nu$; that is, if the measure of $\cyl(\vec a, \bm{a}_i)$ is $\lambda(\cyl(\vec a)) \nu(\bm{a}_i)$. 

On the other hand, if we say that a $p$-tuple $\bm{a}_i$ is \textbf{chosen exceptionally}, then all of the mass associated to $\cyl(\vec a)$ is given to $\cyl (\bm{a}_0, \ldots, \bm{a}_{i-1}, \bm{a}_i)$, where $\bm{a}_i$ is the $p$-tuple $(\rho(\vec a), \rho_2(\vec a), \ldots, \rho_p(\vec a))$.  

We now have the necessary building blocks to describe the construction of the measure $\lambda$. For $1 \leq i < i_1$, the $p$-tuple $\bm{a}_i$ will be chosen typically. Then, the $p$-tuples $\bm{a}_{i_1}, \ldots, \bm{a}_{i_1 + r_1 - 1}$ will be chosen exceptionally. We will define $X_2$ to be the collection of finite sequences $(\bm{a}_0, \ldots, \bm{a}_{i_1 - 1})$ lying in the top half of $\supp \nu^{i_1}$; we will define $X_3$ to be the collection of finite sequences $(\bm{a}_0, \ldots, \bm{a}_{i_1 -1 })$ lying in the bottom half of $\supp \nu^{i_1}$. Notice that $X_2^*$ and $X_3^*$ are disjoint, and each of $X_2^{**}$ and $X_3^{**}$ has been assigned a measure of $\frac{1}{2} + O(2^{-i_1})$.

Next, each $p$-tuple $\bm{a}_i$ for $i_1 + r_1 \leq i < i_2$ will be chosen typically. For those finite sequences $(\bm{a}_0, \ldots, \bm{a}_{i_2 - 1})$ that do not lie in $X_2^*$, choose each $p$-tuple $\bm{a}_{i_2}, \ldots, \bm{a}_{i_2 + r_2 - 1}$ typically. If, instead, $(\bm{a}_0, \ldots, \bm{a}_{i_2 - 1}) \in X_2^*$, choose $\bm{a}_{i_2}, \ldots, \bm{a}_{i_2 + r_2 - 1}$ exceptionally, and let $X_4$ denote the set of those sequences $(\bm{a}_0, \ldots, \bm{a}_{i_2 - 1}) \in X_2^*$ such that $(\bm{a}_{i_1 + r_1}, \ldots, \bm{a}_{i_2 - 1})$ lies in the top half of $\nu^{i_2 - (i_1 + r_1)}$; let $X_5$ denote those sequences $(\bm{a}_0, \ldots, \bm{a}_{i_2 - 1}) \in X_2^*$ such that $(\bm{a}_{i_1 + r_1}, \ldots, \bm{a}_{i_2 - 1})$ lies in the bottom half of $\nu^{i_2 - (i_1 + r_1)}$. Note that $X_4^*$ and $X_5^*$ are disjoint, that $X_4^*$ and $X_5^*$ have each been distributed a mass of $(1/2 \pm O(2^{-i_1})) \cdot (1/2 \pm O(2^{-(i_2 - (i_1 + r_1))})) $, and that $X_4^* \cup X_5^* \subset X_2^*$.

Then, each $p$-tuple $\bm{a}_i$ for $i_2 + r_2 \leq i < i_3$ will be chosen typically. If $(\bm{a}_0, \ldots, \bm{a}_{i_3 - 1}) \notin X_3^*$, then each of $\bm{a}_{i_3}, \ldots, \bm{a}_{i_3 + r_3 - 1}$ will be chosen typically. If, instead, $(\bm{a}_0, \ldots, \bm{a}_{i_3 - 1}) \in X_3^*$, then each of $\bm{a}_{i_3}, \ldots, \bm{a}_{i_3 + r_3 - 1}$ will be chosen exceptionally. In this case, let $X_6$ denote those sequences $(\bm{a}_0, \ldots, \bm{a}_{i_3 - 1}) \in X_3^*$ such that $(\bm{a}_{i_2 + r_2}, \ldots,\bm{a}_{i_3 -1})$ lies in the top half of $\supp \nu^{i_3 - (i_2 + r_2)}$, and let $X_7$ denote those sequences $(\bm{a}_0 \ldots, \bm{a}_{i_3 - 1})$ such that $(\bm{a}_{i_2 + r_2}, \ldots, \bm{a}_{i_3 - 1})$ lies in the bottom half of $\supp \nu^{i_3 - (i_2 + r_2)}$. Note that $X_6^*$ and $X_7^*$ are disjoint, and that $X_6^* \cup X_7^* \subset X_3^*$. Hence $X_4^{**}, X_5^{**}, X_6^{**},$ and $X_7^{**}$ are disjoint sets, each of which has been distributed a mass of $(1/2 \pm O(2^{-i_1})) \cdot (1/2 \pm O(2^{-(i_2 - (i_1 + r_1))}))$.

Now suppose for some $n \geq 4$ that the mass of $\cyl(\vec a)$ has already been chosen for $\vec a$ of the form $\vec a = (\bm{a}_0, \ldots, \bm{a}_{i_{n-1} + r_{n-1}})$, and the sets $X_{1}, \ldots, X_{2n - 1}$ have already been defined. 

We describe the choice of $\bm{a}_i$ for $i_{n-1} + r_{n-1} + 1 \leq i \leq i_n + r_n$ and the choice of $X_{2n}$ and $X_{2n + 1}$.  The tuples $\bm{a}_{i}$ for which $i_{n - 1} + r_{n-1} + 1 \leq i < i_n$ are always chosen typically. If $(\bm{a}_0, \ldots, \bm{a}_{i_n - 1}) \notin X_n^*$, then choose the tuples $\bm{a}_{i_n}$, \ldots, $\bm{a}_{i_n + r_n}$ typically. If $(\bm{a}_0, \ldots, \bm{a}_{i_n - 1}) \in X_n^*$, then choose $\bm{a}_{i_n}, \ldots, \bm{a}_{i_n + r_n}$ exceptionally. In this case, we say that the sequence $(\bm{a}_0, \ldots, \bm{a}_{i_n})$ belongs to $X_{2n}$ if $(\bm{a}_{i_{n - 1} + r_{n-1}}, \ldots, \bm{a}_{i_{n} - 1})$  belongs to the top half of $\supp \nu^{i_{n} - (i_{n - 1} + r_{n - 1})}$; we say that the sequence belongs to $X_{2n + 1}$ if $(\bm{a}_{i_{n - 1} + r_{n - 1}}, \ldots, \bm{a}_{i_n - 1})$ belongs to the bottom half of $\supp \nu^{i_{n} - (i_{n - 1} + r_{n - 1})}$. This means that the total mass of $X_{2n}^{**}$ is equal to
\[\lambda(X_n^{**}) \cdot \left(\frac{1}{2} + O(2^{-(i_n - (i_{n-1} + r_{n-1}))}) \right)\]

Iterating this construction yields a mass distribution $\lambda$ on the cylinder $\sigma$-algebra of infinite sequences of integers. We will write $\supp \lambda$ to refer to the set of infinite sequences that have been chosen according to the above construction. With respect to this mass distribution, each of the disjoint sets $X_{2^t}^{**}, \ldots, X_{2^{t + 1} - 1}^{**}$ has been assigned a measure $2^{-t} + O(2^{-t - i_1})$ with implicit constant independent of $t$; equivalently, the measure of the set $X_n^{**}$ is equal to $w_n(1 + O(2^{-i_1})) \leq 2 w_n$. The construction is arranged so that if $x = (\bm{a}_0, \bm{a}_1, \ldots)$ is an infinite sequence in the support of $\lambda$, then $p$-tuples $\bm{a}_{i_n}, \ldots, \bm{a}_{i_n + r_n}$ will be chosen exceptionally if and only if $x \in X_n^{**}$. Hence, only at most a $2w_n$-fraction of the $\lambda$-mass will be assigned to sequences for which $\bm{a}_{i_n}, \ldots, \bm{a}_{i_n + r_n}$ will be chosen exceptionally. Moreover, each sequence lying in $\supp \lambda$ must belong to infinitely many sets $X_n^{**}$.

We claim that the image of $\supp \lambda$ under $g$ is a closed set, and hence $\supp g^{\sharp} \lambda = g(\supp \lambda)$. Indeed suppose $\{y_j\}_{j=1}^{\infty}$ is a sequence of elements of $g(\supp \lambda)$ that has a limit $y \in \mathbb{R}$. First, it is clear that $y$ must be irrational; if $y$ is rational, then $y$ has a finite continued fraction expansion and cannot be arbitrarily close to numbers with infinite continued fraction expansions.  Moreover, $\{y_j\}_{j=1}^{\infty}$ is a Cauchy sequence; hence, each partial quotient of the continued fraction of $y_j$ must eventually be equal to the corresponding partial quotient of $y$. Thus each finite truncation of the continued fraction expansion of $y$ matches the truncation of an element of $g(\supp \lambda)$, and hence $y \in g(\supp \lambda)$. Hence $g(\supp \lambda)$ is a closed set and $g(\supp \lambda) = \supp g^{\sharp} \lambda$.

Since elements of $\supp \lambda$ belong to infinitely many sets $X_n^*$, it follows that if $y \in \supp g^{\sharp} \lambda$ then there are infinitely many $n$ such that the $p$-tuples $\bm{a}_{i_n}, \ldots, \bm{a}_{i_n + r_n - 1}$ were all chosen exceptionally. This means that $y \in E(S, \infty)$ as desired. So $g^{\sharp} \lambda$ is supported on $E(S, \infty)$.
\section{An estimate on the Fourier transform of $g^{\sharp} \lambda$}
All that remains is to compute an estimate on the Fourier transform of $\lambda$, and it is at this point that a judicious choice of the rapidly growing sequence $\{i_n\}_{n=1}^{\infty}$ and the slowly growing sequence $\{r_n\}_{n=1}^{\infty}$ comes into play. We will estimate $\widehat{g^{\sharp} \lambda}(\xi)$ for large real numbers $\xi$. 

For any $n$, we define
\[\lexcn = \lambda|_{X_{n-1}^{**} \cup X_n^{**} \cup X_{n+1}^{**}}.\]
and
\[\ltypn = \lambda - \lexcn.\]

We will show the following useful fact about continuants for $x \in \ltypn$.

\begin{mylem}\label{ltyp_scale}
Suppose $(\bm{a}_0, \bm{a}_1, \ldots, \bm{a}_i) \notin X_{n-1}^* \cup X_n^*$. If $i_{n} \leq i \leq i_{n + 1}$,
\[\left| \log q(\bm{a}_0, \ldots, \bm{a}_{i}) - i  \sigma \right| \leq \frac{1}{100} i \sigma.\]
Moreover, if $n$ is sufficiently large,
\[\log(\ltypn(\cyl(\bm{a}_0, \ldots, \bm{a}_i))) \leq - \frac{196}{100} i \sigma.\]
\end{mylem}

Let $\alpha = \frac{50}{358}$. Given $\xi$, we choose $i(|\xi|^{\alpha}) = \frac{\log |\xi|^{\alpha}}{\sigma}.$ We fix $n$ such that $i_n \leq i(|\xi|^{\alpha}) \leq i_{n+1}$. Since $n$ will be fixed for the remainder of this argument, we write $\ltyp$ for $\ltypn$ and $\lexc$ for $\lexcn$. 

The total mass $\norm{\lexc}_{\text{TV}} \leq \lambda(X_{n-1}^{**}) + \lambda(X_n^{**}) + \lambda(X_{n+1}^{**}) \leq \frac{6}{n-1}.$
Notice that this decays to zero as $|\xi| \to \infty$. We use only this trivial estimate on $\lexc$.

It remains to estimate $\widehat{g^{\sharp} \ltyp}$. We need the following claim that will be established later.
We will obtain a better estimate for $\widehat{g^{\sharp} \ltyp}$ by mimicking Kaufman's argument. We wish to estimate the Fourier transform
\[\widehat{g^{\sharp} \ltyp}(\xi) = \sum_{x_0}C_{x_0} \int e \left( \xi \frac{p g(x) + p'}{q g(x) + q'} \right) \, d \lambda_{x_0}(x)\]
Here, the sum is extended over finite sequences $x_0 = (\bm{a}_1, \ldots, \bm{a}_{i(|\xi|^{\alpha})})$ such that $\cyl(x_0)$ intersects $\supp \ltyp$, the fractions $\frac{p}{q}$ and $\frac{p'}{q'}$ are the final two convergents to the finite continued fraction $g(x_0)$, the constants $C_{x_0} = \ltyp(\cyl(x_0)) = \lambda(\cyl(x_0))$ sum to $\norm{\ltyp}_{\text{TV}} \leq 1$, and the measures $\lambda_{x_0}$ are the conditional probability measure associated to the sequence $x_0$ defined by $\lambda_{x_0}(A) = \frac{\lambda(x_0 \cdot A)}{\lambda(\cyl(x_0))}.$

We claim that for each $x_0$, we have the $x_0$-independent estimate
\begin{equation}\label{lambdax0lambdax}
\left| \int e \left(\xi \frac{p g(x) + p'}{q g(x) + q'} \right) d(\lambda_{x_0} - \lambda_K)(x) \right| \lesssim |\xi|^{-100}.
\end{equation}
Recall that $x_0 \notin X_{n-1}^* \cup X_n^* \cup X_{n+1}^*$. Hence, if $x = (\bm{a}_{i(|\xi|^{\alpha}) + 1}, \ldots, \bm{a}_{i_{n + 2}}, \ldots) \in \supp \lambda_{x_0}$, then each $p$-tuple of $x$ up to $\bm{a}_{i_{n+2}}$ must have been chosen typically. Since $i(|\xi|^{\alpha}) + 1 \leq i_{n+1}$, it follows that each the first $(i_{n + 2} - i_{n + 1})$ $p$-tuples of $x$ were chosen typically. This means that if $x' = (\bm{a}_{i(|\xi|^{\alpha}) + 1}, \ldots, \bm{a}_{i_{n + 2}})$, then 
\[\lambda_{x_0} (\cyl(x')) = \lambda_K (\cyl(x')).\] 
Because \eqref{qestimate} implies that $q(x') \gtrsim \exp \left(\frac{99}{100} (i_{n +2} - i(|\xi|^{\alpha})) \sigma \right)$, we have that 
\begin{equation}\label{xprime}
\diam(\cyl(x')) \lesssim \exp \left(-\frac{198}{100} (i_{n+2} - i(|\xi|^{\alpha}) \sigma \right) \leq \exp \left( - \frac{198}{100} (i_{n-2} - i_{n-1}) \sigma \right)
\end{equation}

Decomposing the integral in \eqref{lambdax0lambdax} by the value of $x'$, we have the following estimate for sufficiently large $n$:
\begin{IEEEeqnarray*}{Cl}
& \left| \int e \left(\xi \frac{p g(x) + p'}{q g(x) + q'} \right) d(\lambda_{x_0} - \lambda_K)(x) \right| \\
\leq & \sum_{x'} \left| \int_{\cyl(x')} e \left(\xi \frac{p g(x) + p'}{q g(x) + q'} \right) d (\lambda_{x_0} - \lambda_K)(x)\right| \\
\leq & \sum_{x'} \lambda_K(\cyl(x')) |\xi| \sup_{x,y \in \cyl(x')} \left| e\left(\xi \frac{p g(x) + p'}{q g(x) + q'} \right) - e \left(\xi \frac{p g(y) + p'}{q g(y) + q'} \right) \right| \\
\lesssim & \sum_{x'} \lambda_K(\cyl(x')) \diam(\cyl(x')) \exp \left(\frac{\sigma i(\xi)}{\alpha}\right) \sup_{t \in [1, N + 1]} \frac{1}{(q t + q')^2}\\
\lesssim & \sum_{x'} \lambda_K(\cyl(x')) \exp\left (\frac{1}{2} i_{n + 2} \sigma\right)\exp(- \frac{198}{100} (i_{n+2} - i_{n+1}) \sigma)
\lesssim \exp(- i_{n + 2} \sigma),
\end{IEEEeqnarray*}
where above we use the mean value theorem and the fact that the derivative of $\frac{p t + p'}{q t + q'}$ is equal to $\pm \frac{1}{(q t + q')^2},$ which is bounded above in absolute value by $1$.

Hence

\[\widehat{g^{\sharp} \ltyp}(\xi) = \sum_{x_0} C_{x_0} \int e \left(\xi \frac{p g(x) + p'}{q g(x) + q'} \right) \, d \lambda_K(x) + O \left(\exp\left(-i_{n + 2}\sigma \right) \right).\]
But the measure $\lambda_K$ does not depend on $x_0$, so the sum and integral can be interchanged. Recalling that for $n$ sufficiently large we have that $\exp \left(-\frac{1}{4} i_{n+2} \sigma\right) \lesssim |\xi|^{-100}$, we thus have the estimate
\[\widehat{g^{\sharp} \ltyp}(\xi) = \int \sum_{x_0} C_{x_0} e \left( \xi \frac{p t + p'}{q t + q'} \right) \, d g^{\sharp} \lambda_K(t) + O(|\xi|^{-100})\]
We show in Section \ref{KaufmanArgument} the inequality
\begin{equation}\label{gsharplambdatypest}
\left| \int \sum_{x_0} C_{x_0} e \left(\xi \frac{pt + p'}{qt + q'} \right) \, d g^{\sharp} \lambda_K(t) \right| \lesssim |\xi|^{-\frac{1}{100}}.
\end{equation}
Thus
\begin{equation}\label{gsharplambdaest}
|\widehat{g^{\sharp} \lambda}(\xi)| \lesssim \frac{1}{n-1} + |\xi|^{-\epsilon}.
\end{equation}
\section{Proof of Lemma \ref{ltyp_scale}}
\begin{proof}[Proof of Lemma \ref{ltyp_scale}]
First, we prove the estimate for $q(\bm{a}_0, \ldots, \bm{a}_i)$.
We prove the statement by induction on $n$. If $n = 1$, then each $p$-tuple $(\bm{a}_0, \ldots, \bm{a}_i)$ is chosen typically and thus $|\log q(\bm{a}_0, \ldots, \bm{a}_i) - \sigma i| \leq \frac{1}{10000} \sigma i$, as desired.
Suppose instead $n \geq 2$. Suppose $x \in \supp \ltypn$. Let 
\[x' = (\bm{a}_0, \ldots, \bm{a}_{i_{n-1}}, \ldots, \bm{a}_i)\]
be the initial sequence of $x$. Because $x \notin X_{n-1}^{**} \cup X_n^{**}$, we know that each $p$-tuple after $\bm{a}_{i_{n-2} + r_{n-2} - 1}$ is chosen typically. By Lemma \ref{joining}, it follows that 
\[|\log q(\bm{a}_0, \ldots, \bm{a}_{i}) - \log q(\bm{a}_{i_{n-2} + r_{n-2}}, \ldots, \bm{a}_{i}) |  \leq  C_N + |\log q (\bm{a}_0, \ldots, \bm{a}_{i_{n-2} + r_{n-2} - 1})|.\]
Since each $p$-tuple in $(\bm{a}_{i_{n-2}+ r_{n-2}}, \ldots, \bm{a}_i)$ is chosen typically, we have from \eqref{qestimate} that
\[\left| \log q(\bm{a}_{i_{n-2} + r_{n-2}}, \ldots, \bm{a}_i) - \sigma (i - (i_{n-2} + r_{n-2}))\right| \leq \frac{1}{1000} \sigma (i - (i_{n-2} + r_{n-2})) \]
Since $i - (i_{n-2} + r_{n-2}) \geq i_n - i_{n-1} \geq 100i_{n-2}$ for sufficiently large $n$, we have that  
\[\left| \log q(\bm{a}_{i_{n-2} + r_{n-2}}, \ldots, \bm{a}_{i}) - \sigma i \right| \leq \frac{1}{100} \sigma i\]

It remains to estimate the logarithm of the denominator of the initial sequence
\[\log q(\bm{a}_0, \ldots, \bm{a}_{i_{n-2} + r_{n-2} - 1}).\]
If $x \notin X_{n-2}^{**}$, then we can apply the induction hypothesis to conclude that 
\[\log q(\bm{a}_0, \ldots, \bm{a}_{i_{n-2} + r_{n-2} - 1}) \leq 2 \sigma i_{n-2}.\]

If $x \in X_{n-2}^{**}$, we recall that each of the $p$-tuples $\bm{a}_{i_{n-2}}, \ldots, \bm{a}_{i_{n-2} + r_{n-2} - 1}$ is chosen exceptionally, so 
\[q(\bm{a}_0, \ldots, \bm{a}_{i_{n-2}}, \bm{a}_{i_{n-2} + r_{n-2} - 1}) \leq \phi_{pr_{n-2} }(\bm{a}_0, \ldots, \bm{a}_{i_{n-2}})\]
Upon taking logarithms, we conclude that
\[\log q(\bm{a}_0, \ldots, \bm{a}_{i_{n-2}}, \bm{a}_{i_{n-2} + r_{n-2} - 1}) \leq \log \phi_{p r_{n-2}} ((\bm{a}_0, \ldots, \bm{a}_{i_{n-2}})),\]
and the conditions on the $i_n$ imply that this is is no more than
\[\frac{1}{1000} \sigma (i_{n} - i_{n-1}) \leq \frac{1}{1000} \sigma i\]
as desired.

It remains to prove the estimate on $\lambda(\cyl(\bm{a}_0, \ldots, \bm{a}_i))$. By the definition of $\lambda$ and \eqref{nuestimate}, we have
\begin{IEEEeqnarray*}{rCl}
\lambda(\cyl(\bm{a}_0, \ldots, \bm{a}_i)) & = & \lambda(\cyl(\bm{a}_0, \ldots, \bm{a}_{i_{n-2} + r_{n-2} - 1})) \nu(\bm{a}_{i_{n-2} + r_{n-2}}) \cdots \nu(\bm{a}_i) \\
& \leq & \nu(\bm{a}_{i_{n-2} + r_{n-2}}) \cdots \nu(\bm{a}_i) \\
& \leq & \exp \left(- \frac{198}{100} \sigma (i - (i_{n-2} + r_{n-2})) \right) 
\end{IEEEeqnarray*}
If $n$ is sufficiently large, we have $i \geq i_n \geq 1000 i_{n-1} \geq 1000 i_{n-2 + r_{n-2}}$, so $i - (i_{n-2} + r_{n-2}) \geq \frac{999}{1000} i$. Thus 
\[\lambda(\cyl(\bm{a}_0, \ldots, \bm{a}_i)) \leq \exp \left(-\frac{196}{100} \sigma i \right),\]
as desired.
\end{proof}
\section{Proof of the estimate for $\widehat{g^{\sharp} \ltyp}$}\label{KaufmanArgument}
\subsection{Oscillatory integral estimates of Kaufman}
Our proof of the estimate \eqref{gsharplambdatypest} follows Kaufman \cite{Kaufman80}. Kaufman's argument is based on a few integral inequalities. We summarize these below. All three of these inequalities appear in Kaufman's article \cite{Kaufman80}.

\begin{mylem}[Integral inequality from \cite{Kaufman80}]\label{intineq}
Let $f$ be a function such that $|f(t)| \leq 1$ and $|f'(t)| \leq M$ on an interval $[a,b]$, and write $m_2 = \int_a^b |f(t)|^2 \, dt$. Let $\lambda_K$ be a probability measure on $[a,b]$ and let $\Omega(u)$ be the maximum $\lambda_K$-measure of all intervals $[t, t + u] \subset [a,b]$ for $u > 0$. Then
\[\int_a^b |f(t)| \, d \lambda_K(t) \leq 2 M^{1/10} m_2^{3/10} + \Omega(M^{-9/10} m_2^{3/10}) (1 + M^{7/10} m_2^{1/10}).\]
\end{mylem}

In practice, applying this inequality requires two van der Corput type lemmas to estimate $m_2$. The first of these lemmas is a non-stationary phase estimate that will be useful for tuples $\vec a_1, \vec a_2$ such that the second-last denominators $q'(\vec a_1)$ and $q'(\vec a_2)$ agree.

\begin{mylem}\label{nonstationary}[Non-stationary oscillatory integral estimate, \cite{Kaufman80}]
Let $f \in C^2[0,1]$, and suppose $f' \geq a > 0$ or $f' \leq -a < 0$, and $|f''| < b$ on the interval $[0,1]$. Then
\[\left| \int e(h(x)) \, dx \right| < a^{-1} + a^{-2} b.\]
\end{mylem}

The second lemma is a stationary phase estimate that is useful for pairs of consecutive denominators where $q'(\vec a_1)$ and $q'(\vec a_2)$ differ.

\begin{mylem}\label{stationary}[Stationary oscillatory integral estimate, \cite{Kaufman80}]
Let $h \in C^2[0,1]$ be such that $h'(x) = (a_1 x + a_2)g(x)$ on $[0,1]$, where $g \in C^1$, $|g(x)| \geq a$, $|g'(x)| \leq b$ on $[0,1]$, where $b > 1$. Then
\[\left| \int_0^1 e(h(x)) \, dx \right| < 6 b a^{-3/2} |a_1|^{-1/2}.\]
\end{mylem}
\subsection{Estimate of integral}
We seek to estimate the integral
\[\int f(t) \, d \lambda_K(t) ,\]
where 
\[f(t) = \sum_{x_0 : \cyl(x_0) \cap \supp \ltyp \neq \emptyset} \lambda(\cyl (x_0)) e \left(\xi \frac{p(x_0) t + p'(x_0)}{q(x_0) t + q'(x_0)} \right).\]
Because the sum over $x_0$ is extended over those $x_0$ such that $\cyl(x_0) \cap \supp \ltyp \neq \emptyset$, it follows that $x_0 \notin X_{n-1}^* \cup X_n^* \cup X_{n+1}^*$. Therefore, by Lemma \ref{ltyp_scale}, it follows for each $x_0$ in the sum that
\begin{equation}\label{qx0estimate}
|\xi|^{\frac{99}{100} \alpha} \lesssim q'(x_0), q(x_0) \lesssim |\xi|^{\frac{101}{100} \alpha}.
\end{equation}

\begin{IEEEeqnarray}{rCl}
M & \sim & |\xi|^{\frac{244}{358}} \label{Mbound}\\
m_2 & \lesssim & |\xi|^{-\frac{99}{358}} \label{m2bound}\\
\Omega(u) & \lesssim & u^{98/100} \label{Omegabound},
\end{IEEEeqnarray}
giving $\int f(t) d \lambda_K(t) \lesssim |\xi|^{-\frac{1}{100}}.$

The estimate on $\Omega(u)$ is the estimate \eqref{lambdafrostman}. We outline the estimate for $M$ and $m_2$ below.
\subsection{Estimate for $M$}
First, we must obtain an upper bound on $|f'(t)|$. Because $\left|p(x_0) q'(x_0) - p'(x_0) q(x_0)\right| = 1$, we have from the triangle inequality that
\[\left|f'(t) \right| \leq |\xi| \sum_{x_0} \lambda(\cyl(x_0)) \frac{1}{(q(x_0) t + q'(x_0))^2}\]
By \eqref{qx0estimate},  we conclude
\[\left|f'(t) \right| \lesssim |\xi|^{1 - \frac{198}{100} \alpha}.\]
Plugging in $\alpha = \frac{50}{358}$ gives \eqref{Mbound}.
\subsection{Estimate of $L^2$ norm}
It remains to estimate $m_2$. By expanding the $L^2$-norm, we see that
\begin{equation}\label{m2sum}
m_2 = \sum_{x_0, y_0} \lambda(\cyl(x_0)) \lambda(\cyl(y_0)) \int_{1}^{N+1} e \left( \xi \left( \frac{p(x_0) t + p'(x_0)}{q(x_0) t + q'(x_0)} - \frac{p(y_0)t + p'(y_0)}{q(y_0) t + q_0'(y_0)} \right) \right) \, dt
\end{equation}
The derivative of the argument with respect to $t$ can be written in the form 
\[\xi \frac{((q(y_0) + q(x_0))t + q'(y_0) + q'(x_0))}{(q(x_0) + t q'(x_0))^2(q(y_0) + t q'(y_0))^2}  ((q(y_0) - q(x_0))t + (q'(y_0) - q'(x_0))).\]
Let $l(t)$ denote the cofactor
\[l(t) = \xi  \frac{((q(y_0) + q(x_0))t + q'(y_0) + q'(x_0))}{(q(x_0) + t q'(x_0))^2(q(y_0) + t q'(y_0))^2}.\]

A straightforward calculation using the estimate \eqref{qx0estimate} shows that 
\[|l(t)| \gtrsim |\xi|^{1 - \frac{304}{100} \alpha}; \qquad |l'(t)| \lesssim |\xi|^{1 - \frac{292}{100} \alpha}.\]

We split the sum over $(x_0, y_0)$ into three different sets: $S_1$, the set of pairs $(x_0, y_0)$ such that $q(x_0) = q(y_0)$ but $q'(x_0) \neq q'(y_0)$, $S_2$, the set of pairs $(x_0, y_0)$ such that $q(x_0) \neq q(y_0)$, and $S_3$, the set of pairs $(x_0,y_0)$ such that $q(x_0) = q(y_0)$ and $q'(x_0) = q'(y_0)$. We write $\Sigma_j$ for the sum over $S_j$.

We will show the estimates 
\begin{IEEEeqnarray}{rCl}
\Sigma_1 & \lesssim & |\xi|^{-(1 - \frac{316}{100} \alpha)} \label{Sigma1bound} \\
\Sigma_2 & \lesssim & |\xi|^{-\frac{1}{2} (1 - \frac{328}{100} \alpha)} \label{Sigma2bound} \\
\Sigma_3 & \lesssim & |\xi|^{-\frac{194}{100} \alpha}. \label{Sigma3bound}
\end{IEEEeqnarray}
Combining the estimates \eqref{Sigma1bound}, \eqref{Sigma2bound}, and \eqref{Sigma3bound} together with the choice $\alpha = \frac{50}{358}$ gives the estimate $m_2 \lesssim |\xi|^{-\frac{97}{358}}$.
\subsubsection{Estimate of $\Sigma_1$}
We begin with the sum $\Sigma_1$. We must estimate
\[\sum_{(x_0, y_0) \in S_1} \lambda(\cyl(x_0)) \lambda(\cyl(y_0)) \int_1^{N_1} e \left( \xi \left( \frac{p(x_0) t + p'(x_0)}{q(x_0) t + q'(x_0)} - \frac{p(y_0)t + p'(y_0)}{q(y_0) t + q_0'(y_0)} \right) \right) \, dt.\]
Because $(x_0, y_0) \in S_1$, we have that $q(x_0) = q(y_0)$ but $q'(x_0) \neq q'(y_0)$; hence, the derivative of the phase is given by
\[l(t) (q'(y_0) - q'(x_0))\]
where $|q'(y_0) - q'(x_0)| \geq 1$. Hence we can apply Lemma \ref{nonstationary} to estimate each integral, where $a \gtrsim |\xi|^{1 - \frac{304}{100} \alpha}$ and $b \lesssim |\xi|^{1 - \frac{292}{100} \alpha}$. Since $\sum_{(x_0, y_0) \in S_1} \lambda(\cyl(x_0)) \lambda(\cyl(y_0)) \leq 1$, we have 
\[\Sigma_1 \lesssim |\xi|^{-(1 - \frac{316}{100} \alpha)}.\]
This establishes the bound \eqref{Sigma1bound}.
\subsubsection{Estimate of $\Sigma_2$}
Next we estimate the sum $\Sigma_2$, given by
\[\sum_{(x_0, y_0) \in S_2} \lambda(\cyl(x_0)) \lambda(\cyl(y_0)) \int_1^{N+1} e \left( \xi \left( \frac{p(x_0) t + p'(x_0)}{q(x_0) t + q'(x_0)} - \frac{p(y_0)t + p'(y_0)}{q(y_0) t + q_0'(y_0)} \right) \right) \, dt.\]
This time, because $(x_0, y_0) \in S_2$, the derivative of the phase is given by 
\[l(t) ((q(y_0) - q(x_0))t + q'(y_0) - q'(x_0)).\]
Because the coefficient $q(y_0) - q(x_0) \neq 0$, we apply Lemma \ref{stationary} to estimate each integral; we use the fact that $|a_1| = |q(y_0) - q(x_0)| \geq 1$, that $a \gtrsim |\xi|^{1 - \frac{304}{100} \alpha}$, and that $b \lesssim |\xi|^{1 - \frac{292}{100} \alpha}$. Thus Lemma \ref{stationary} gives that each integral in the expression for $\Sigma_2$ is bounded above by $\lesssim |\xi|^{-\frac{1}{2} (1 - \frac{328}{100} \alpha)}$. Because $\sum_{(x_0, y_0) \in S_2} \lambda(\cyl(x_0)) \lambda(\cyl(y_0)) \leq 1$, we have 
\[\Sigma_2 \lesssim |\xi|^{-\frac{1}{2}(1 -\frac{328}{100} \alpha)}.\]
This shows the bound \eqref{Sigma2bound}.
\subsubsection{Estimate of $\Sigma_3$}
Finally, we estimate $\Sigma_3$:
\[\sum_{(x_0, y_0) \in S_3} \lambda(\cyl(x_0)) \lambda(\cyl(y_0)) \int_1^{N+1} e \left( \xi \left( \frac{p(x_0) t + p'(x_0)}{q(x_0) t + q'(x_0)} - \frac{p(y_0)t + p'(y_0)}{q(y_0) t + q_0'(y_0)} \right) \right) \, dt.\]
For this sum, we simply estimate each integral above by $N \lesssim 1$. Now, $S_3$ consists of those pairs $(x_0, y_0)$ such that $q(x_0) = q(y_0)$ and $q'(x_0) = q'(y_0)$. Because $x_0$ is, up to the integer part, determined by $q(x_0)$ and $q'(x_0)$, the sum is controlled (up to a constant depending only on $N$) by the diagonal terms $x_0 = y_0$. Hence it is enough to estimate
\[\sum_{x_0: \cyl(x_0) \cap \supp \ltyp \neq \emptyset} \lambda(\cyl(x_0))^2.\]
For such $x_0$, we have estimated in Lemma \ref{ltyp_scale} that $\lambda(\cyl(x_0)) \lesssim \exp \left(-\frac{196}{100} \sigma i(|\xi|^{\alpha}) \right) \lesssim |\xi|^{-\frac{194}{100}\alpha}$. Hence
\begin{IEEEeqnarray*}{Cl}
& \sum_{x_0: \cyl(x_0) \cap \supp \ltyp \neq \emptyset} \lambda(\cyl(x_0))^2 \\
\lesssim & |\xi|^{-\frac{194}{100} \alpha} \sum_{x_0 : \cyl(x_0) \cap \supp \ltyp \neq \emptyset} \lambda(\cyl(x_0)) \\
\leq & |\xi|^{-\frac{194}{100} \alpha}.
\end{IEEEeqnarray*}
This proves the bound \eqref{Sigma3bound}, completing the estimate for $m_2$. 
\section{Quantitative statements for a few specific examples}
The conditions (A)-(D) are the only conditions required on $i_n$ and $r_n$ for the argument to work. This allows us to turn \eqref{gsharplambdaest} into a quantitative statement. Lemma \ref{ltyp_scale} allows us to obtain useful quantitative bounds on the growth rate of $i_n$ required for the argument to work.

We will specialize to the case $E(\psi, \infty)$ for the rest of this section.

Let $\psi(q)$ be a function such that $q^2 \psi(q)$ is decreasing to zero. We will choose the $i_n$ to satisfy the following growth condition:
\begin{equation}\label{Epsiinftycond}
i_{n+1} \geq \log (\Phi^{p r_n} (2^{p r_n - 1} \exp (2 \sigma i_n)))
\end{equation}
where $\Phi(q) = \frac{1}{q \psi(q)}$.
Because $\frac{\Phi(q)}{q} \to \infty$, The inequality \eqref{Epsiinftycond} implies (A).

To check (B), we observe that
\begin{IEEEeqnarray*}{rCll}
i_{n+1}& \geq & \log \left( \Phi^{p r_n} (2^{p r_n - 1} \exp(2 \sigma i_n)) \right) \quad & \\
& \geq & \log \Phi^{p r_n} (2^{p r_n - 1} q(\vec a)) \quad & \text{for $\vec a = (\bm{a}_0, \ldots, \bm{a}_{i_n}) \in X_n^*$ by Lemma \ref{ltyp_scale}.}\\
& \geq & \log \phi_{p r_n}(\vec a) \quad & \text{for $\vec a = (\bm{a}_0, \ldots, \bm{a}_{i_n})\in X_n^*$ by \eqref{phirestimate}.}
\end{IEEEeqnarray*}
Using (A), we know that $i_{n + 2} - i_{n + 1} \geq i_{n+1}$ for sufficiently large $n$, establishing (B). Hence the condition \eqref{Epsiinftycond} implies (A) and (B).

We now give a few interesting special cases. 
\begin{myex} 
If $\psi(q) = q^{-\tau}$ for some $\tau > 2$, and $\omega$ is any function that increases to $\infty$, then we can choose $\lambda$ so that
\[\left|\widehat{g^{\sharp} \lambda} (\xi) \right| \lesssim \frac{\omega(|\xi|)}{\log \log \xi}.\]
\end{myex} 
\begin{proof}
If $\psi(q) = q^{-\tau}$, then $\Phi(q) = q^{\tau - 1}$, so the condition \eqref{Epsiinftycond} becomes
\[i_{n + 1} \geq (\tau - 1)^{p r_n}((p r_n- 1) \log 2 + 2 \sigma i_n).\]
If $i_n \gg r_n$, then this holds if we choose
\[i_{n +1} = \lceil 3 (\tau - 1)^{pr_n} \sigma i_n \rceil.\] 
If $\tilde \omega(n)$ is any function that increases to $\infty$, we can choose $r_n$ to increase sufficiently slowly depending on $\tilde \omega(n)$ so that for sufficiently large $n$, we have
\[i_n \leq \exp(n \tilde \omega(n)).\]
Solving for $n$, we see
\[n \tilde \omega(n) \geq \log i_n.\]
But recall that for a given $\xi$, we have 
\[i(|\xi|^{\alpha}) \sim \log |\xi|\]
so if $i_n < i(|\xi|^{\alpha}) \leq i_{n + 1}$, then
\[n \tilde \omega(n) \geq \log \log |\xi|.\]
This is the same thing as saying that, if $\omega$ is any function that increases to $\infty$, we can choose $i_n$ such that
\[n \geq \frac{\log \log |\xi|}{\omega(|\xi|)}.\]
Hence by \eqref{gsharplambdaest}, we chose
\[|\widehat{g^{\sharp} \lambda}(\xi)| \lesssim \frac{\omega(\xi)}{\log \log |\xi|}.\]
\end{proof}
This falls short of the condition \eqref{normalcond} that guarantees the existence of normal numbers in $E(\psi, \infty)$ in two ways: not only is the power on $\log \log |\xi|$ equal to $1$, but there is an additional loss of a slowly growing function in the numerator.

\begin{myex}
If $\psi(q) = \exp(-q)$, and $\omega$ is any function that increases to $\infty$, we can choose $\lambda$ so that
\[\left|\widehat{g^{\sharp}\lambda}(\xi)\right| \lesssim  \frac{\omega(|\xi|)}{\log^*(|\xi|)},\]
where $\log^*(|\xi|)$ denotes the minimal integer $t$ such that $0 < \log^t(|\xi|) \leq 1$.
\end{myex}
\begin{proof}
This time, $\Phi(q) = \exp(q)/q \leq \exp(q)$, so the conditions (A) and (B) will be implied by the condition
\[i_{n + 1} \geq \log \exp^{p r_n} (2^{p r_n - 1}  \exp(\sigma i_n)).\]
This will hold if 
\[i_{n + 1} \geq \exp^{p r_n} (2^{p r_n} \sigma i_n).\]
If the $r_n$ grow slowly enough relative to the $i_n$, this is implied by the condition 
\[i_{n +1} \geq \exp^{p r_n + 1} (\sigma i_n).\]
This holds provided that
\[\log^{\sum_{j=1}^n p r_j + 1} i_{n+1} \gtrsim  i_1 \sim 1.\]
Rearranging, we see that it suffices to select the $i_n$ such that
\[\sum_{j=1}^n p r_j + 1 \lesssim \log^*{i_{n}}\]
Because the $r_j$ can increase arbitrarily slowly, we see that if $\tilde \omega$ is any function that increases to $\infty$, it is possible to choose the $i_n$ so that, for sufficiently large $n$,
\[n \tilde \omega(n) \gtrsim \log^* i_n.\]
Since $\log |\xi| \sim i(\xi)$, we have that if $i_n < i(|\xi|^{\alpha}) \leq i_{n +1}$, then it it is possible to choose $\lambda$ such that
\[n \tilde \omega(n) \gtrsim \log^*(\xi).\]
This is equivalent to saying that if $\omega$ is any function that increases to $\infty$, then it is possible to choose the $i_n$ so that
\[n \gtrsim \frac{\log^* |\xi|}{\omega(|\xi|)}.\]
Hence by \eqref{gsharplambdaest}, we have for this choice of $\{r_n\}$ and $\{i_n\}$ that
\[\left| \widehat{g^{\sharp} \lambda}(\xi) \right| \lesssim \frac{\omega(|\xi|)}{\log^*(|\xi|)}.\]
\end{proof}
\bibliographystyle{plain}
\bibliography{Master_Bibliography}
\end{document}